\documentclass[12pt,a4]{amsart}
\usepackage{amsmath}
\usepackage{amssymb}
\usepackage{amsthm}
\usepackage{enumitem}
\usepackage{url}
\usepackage{soul}

\usepackage{color}
\usepackage{dsfont}
\usepackage{epsfig}
\usepackage[hidelinks]{hyperref}
\usepackage{setspace}
\usepackage[authoryear,round,sort]{natbib}
\usepackage{comment}
\usepackage{booktabs}
\usepackage{float}
\usepackage{stmaryrd}
\numberwithin{equation}{section}

\usepackage{geometry}
\theoremstyle{plain}
\newtheorem{proposition}{Proposition}[section]

\newtheorem{lemma}{Lemma}[section]
\newtheorem{corollary}{Corollary}[section]

\theoremstyle{definition}
\newtheorem{definition}{Definition}[section]

\theoremstyle{remark}
\newtheorem{rk}{Remark}[section]
\expandafter\let\expandafter\oldproof\csname\string\proof\endcsname
\let\oldendproof\endproof
\renewenvironment{proof}[1][\proofname]{%
  \oldproof[\noindent\textbf{#1.} ]%
}{\oldendproof}
\newcommand{\1}{\mathds{1}}

\newcommand{\E}{\mathbb{E}}

\newcommand{\be}{\begin{equation}}
\newcommand{\ee}{\end{equation}}
\newcommand{\by}{\begin{eqnarray*}}
\newcommand{\ey}{\end{eqnarray*}}

\renewcommand{\leq}{\leqslant}
\renewcommand{\geq}{\geqslant}
\usepackage{xcolor}
\definecolor{dark-red}{rgb}{0.4,0.15,0.15}
\definecolor{dark-blue}{rgb}{0.15,0.15,0.4}
\definecolor{medium-blue}{rgb}{0,0,0.5}
\hypersetup{
    colorlinks, linkcolor={blue},
    citecolor={blue}, urlcolor={blue}
}
\allowdisplaybreaks

\begin{document}
	\title{On resistance distance of Markov chain and its sum rules}
	\author{Michael C.H. Choi}
	\address{Institute for Data and Decision Analytics, The Chinese University of Hong Kong, Shenzhen, Guangdong, 518172, P.R. China}
	\email{michaelchoi@cuhk.edu.cn}
	\date{\today}
	\maketitle
	
	\begin{abstract}
		Motivated by the notion of resistance distance on graph, we define a new resistance distance between two states on a given finite ergodic Markov chain based on its fundamental matrix. We prove a few equivalent formulations and discuss its relation with other parameters of the Markov chain such as its group inverse, stationary distribution, eigenvalues or hitting time. In addition, building upon existing sum rules for the hitting time of Markov chain, we give sum rules of this new resistance distance of Markov chains that resembles the sum rules of the resistance distance on graph. This yields Markov chain counterparts of various classical formulae such as Foster's first formula or the Kirchhoff index formulae.
		\smallskip
		
		
	\end{abstract}
	
	
	
	\section{Introduction and main results}
	
	On a simple connected graph $G = (V,E)$, the resistance distance $\Omega_{i,j}^G$ between two vertices $i,j \in V$ is defined to be the voltage when a unit current enters $i$ and leaves $j$, see e.g. \cite{Tetali91}. Equivalently, it can be defined via the notion of the generalized inverse $L^{\#} = (L^{\#}_{i,j})_{i,j \in V}$ of the Laplacian $L := D - A$, where $D$ is the diagonal matrix of vertex degrees, $A$ is the adjacency matrix of $G$ and $LL^{\#}L = L$. More precisely, according to \cite{K02,B10}, we have
	\begin{align}\label{eq:omegag}
	\Omega_{i,j}^G := L^{\#}_{i,i} + L^{\#}_{j,j} - L^{\#}_{i,j} - L^{\#}_{j,i}.
	\end{align}
	Motivated by this definition of resistance distance on graph, we would like to define an analogous notion of resistance distance that would play a similar role between two states of a discrete-time homogeneous finite Markov chain $X = (X_n)_{n \in \mathbb{N}_0}$, where we denote $\mathbb{N}_0$ to be the set of non-negative integers. Throughout this article, we consider an ergodic (i.e. irreducible and aperiodic) Markov chain $X$ on a finite state space $\mathcal{X}$ with transition matrix $P = (P_{i,j})_{i,j \in \mathcal{X}}$ and stationary distribution $\pi = (\pi_i)_{i \in \mathcal{X}}$, which is considered to be a row vector of size $|\mathcal{X}|$. Writing $\Pi$ to be the matrix where each row is $\pi$,  the fundamental matrix $F = (F_{i,j})_{i,j \in \mathcal{X}}$ associated with the Markov chain $X$, first proposed in the work of \cite{KS76}, is given by
	$$F := (I - P + \Pi)^{-1},$$
	where $I$ is the identity matrix of size $|\mathcal{X}|\times|\mathcal{X}|$.	Note that the above inverse always exists. In view of \eqref{eq:omegag}, we now define a new resistance distance $\Omega = (\Omega_{i,j})_{i,j\in \mathcal{X}}$ of Markov chain by simply replacing $L^{\#}$ by $F$, that is,
	\begin{definition}[Resistance distance of Markov chain]\label{def:omega}
		Given an ergodic Markov chain $X$ with fundamental matrix $F = (F_{i,j})_{i,j \in \mathcal{X}}$, we define the resistance distance $\Omega_{i,j}$ between two states $i, j \in \mathcal{X}$ to be
		$$\Omega_{i,j} := F_{i,i} + F_{j,j} - F_{i,j} - F_{j,i}.$$
	\end{definition}
	It turns out that this definition of resistance distance admits a few equivalent formulations in terms of other important quantities and parameters of Markov chain, such as the group inverse of $I - P$ as well as the mean hitting time of $X$, see Proposition \ref{prop:equivalentresis} below. To this end, let us proceed by briefly recalling these notions. The group inverse $D = (D_{i,j})_{i,j \in \mathcal{X}}$ of $I - P$, first studied by \cite{M75} in a Markov chain setting, is defined to be the matrix that satisfies
	$$(I-P)D(I-P) = I - P, \quad D(I-P)D = D, \quad (I-P)D = D(I-P).$$
	In this paper, as we only discuss the case where the Markov chain $X$ is ergodic, $D$ can be conveniently expressed as
	$$D = \sum_{n \geq 0} (P^n - \Pi),$$
	see e.g. \cite[Theorem $2.4$]{M75}. The group inverse $D$ also appears under different names in the literature, ranging from deviation matrix \cite{CvD02}, ergodic potential \cite{Syski78} to centered resolvent \cite{Miclo16}. We remark that the notion of group inverse is first introduced in the work of \cite{Erdelyi67}, and group inverse is the special case of Drazin inverse when the index of the matrix is either $1$ or $0$. We now move on to discuss a few probabilistic parameters of interest. For $j \in \mathcal{X}$, we write $\tau_j := \inf\{n \geq 0;~X_n = j\}$ to be the first time that the Markov chain $X$ hits the state $j$, and the usual convention of $\inf \emptyset = \infty$ applies. We also denote $\mathbb{E}_i$ to be the expectation under $X_0 = i$. For example, $\E_i(\tau_j)$ is the mean hitting time of $j$ starting from $i$. Among various hitting time parameters as studied in \cite{AF14}, we are interested in the following three:
	\begin{itemize}
		\item Commute time $t^c_{i,j}$ between $i$ and $j$:
		$$t^c_{i,j} := \E_i(\tau_j) + \E_j(\tau_i).$$
		Note that the commute time defines a metric on $\mathcal{X}$.
		\item Average hitting time $t^{av}$:
		$$t^{av} := \sum_{i,j \in \mathcal{X}} \E_i(\tau_j) \pi_i \pi_j.$$
		$t^{av}$ represents the average hitting time from $i$ to $j$ of $X$ when we sample these two states $i,j$ independently from $\pi$. Note that the average hitting time is equal to Kemeny's constant, that is,
		$$t^{av} = \sum_{j \in \mathcal{X}} \E_i(\tau_j) \pi_j,$$
		where the right hand side is the Kemeny's constant which is independent of the starting state $i$. We refer interested readers to \cite{LL02,PT18,Mao04,Kirk14,CuiMao10} for further references on this parameter.
		\item Forest representation of mean hitting time:
		Let $G(P)$ be the weighted direct graph on vertices $\mathcal{X}$ and arc weights to be the corresponding transition probabilities. The weight of a weighted direct graph is the product of its arc weights, and the weight of a set of weighted direct graphs is the sum of the weights of its members. Define $f_{i,j}$ to be the total weight of $2$-tree in-forests of $G(P)$ that have one tree containing $i$ and the other rooted at $j$, where we recall an in-forest is a spanning subdigraph of $G(P)$ all of whose weak components are converging trees (also known as in-arborescences). Let $q_j$ to be total weight of in-trees rooted at $j$ and $q := \sum_{j \in \mathcal{X}} q_j$. According to the Markov chain tree theorem \cite{AT89} and recent results in \cite{C07,CD18}, one can express the stationary distribution and mean hitting time via these graph-theoretic parameters as, for $i,j \in \mathcal{X}$,
		\begin{align}\label{eq:tree}
			\pi_j &= \dfrac{q_j}{q}, \quad \E_i(\tau_j) = \dfrac{f_{i,j}}{q_j}.
		\end{align} 
	\end{itemize}
	
	With the above notations in mind, we are now ready to present our first result that gives a few equivalent formulations of $\Omega_{i,j}$. These formulations are particularly useful when it comes to proving various properties and sum rules of $\Omega$.
	
	\begin{proposition}\label{prop:equivalentresis}
		The resistance distance $\Omega_{i,j}$ of $X$, as defined in Definition \ref{def:omega}, can be written as, for $i,j \in \mathcal{X}$,
		\begin{enumerate}
			\item(Group inverse representation)\label{it:gir}
			$$\Omega_{i,j} = D_{i,i} + D_{j,j} - D_{i,j} - D_{j,i}.$$
			\item(Mean hitting time representation)\label{it:mhtr}
			$$\Omega_{i,j} = \pi_j \mathbb{E}_i(\tau_j) + \pi_i \mathbb{E}_j(\tau_i).$$
			\item(Forest representation)\label{it:forest}
			$$\Omega_{i,j} = \dfrac{f_{i,j}+f_{j,i}}{q}.$$
			\item(Commute time representation for doubly stochastic $P$)\label{it:ds} 
			When $P$ is doubly stochastic, that is both the row sums and column sums of $P$ are $1$, then we have
			$$\Omega_{i,j} = \dfrac{1}{|\mathcal{X}|}t^{c}_{i,j}.$$			
			In other words, the resistance distance is a scaled version of commute time in the doubly stochastic case.
		\end{enumerate}
	\end{proposition}

	\begin{rk}[Connections with existing notions of resistance distance on weighted direct graph]\label{rk:connections}
		In this Remark, we would like to point out to readers on possible connections with existing notions of resistance distance on weighted direct graph. As $L(P) := I - P$ can be interpreted as the Laplacian matrix of the weighted direct graph corresponding to the Markov chain $X$ (see e.g. \cite[Section $2.2$]{CA02}), existing notions of effective resistance on directed graph are thus closely related to the proposed resistance distance $\Omega$.
		\newline
		In \cite{YSL16a, YSL16b}, the authors propose a notion of effective resistance $R = (R_{i,j})_{i,j \in \mathcal{X}}$ on weighted direct graph via the reduced Laplacian. Precisely, let $\1_N$ be the all-ones vector of length $N := |\mathcal{X}|$ and let $I_{N}$ be the identity matrix of size $N$. Let $Q \in \mathbb{R}^{(N-1)\times N}$ be any matrix that satisfies
		$$Q\1_N = 0, \quad QQ^T = I_{N-1}, \quad Q^TQ = I_N - \dfrac{1}{N} \1_N \1_N^T.$$
		Reduced Laplacian $\overline{L(P)}$ of $L(P)$ is then defined to be
		$$\overline{L(P)} := Q L(P) Q^T.$$
		Let $\Sigma$ be the unique solution to the Lyapunov equation
		$$\overline{L(P)} \Sigma + \Sigma \overline{L(P)}^T = I_{N-1},$$
		and $X$ be
		$$Y := 2 Q^T \Sigma Q.$$
		The authors in \cite{YSL16a,YSL16b} define the effective resistance to be
		$$R_{i,j} := Y_{i,i} + Y_{j,j} - 2Y_{i,j}.$$
		Note that according to \cite[Section II]{YSL16a} $\overline{L(P)}$ is not unique and depends on the choice of $Q$, while $R$ is independent of the choice of $Q$. To compare $R$ with our proposed resistance distance $\Omega$, it boils down to a comparison between $Y$ and the fundamental matrix $F$ or the group inverse $D$. While both $R$ and $\Omega$ does not define a metric in general, in Proposition \ref{prop:rmetric} below we show that $\Omega$ does define a metric when $P$ is doubly stochastic while it is unclear whether $R$ also defines a metric in this setting. On the other hand however, $\sqrt{R} = (\sqrt{R_{i,j}})$ defines a metric (see \cite[Theorem $3$]{YSL16a}) while it is not clear whether $\sqrt{\Omega} = (\sqrt{\Omega_{i,j}})$ defines a metric. In Section $V$ of \cite{YSL16a}, the authors motivate their definition $R$ by outlining a few drawbacks in defining resistance distance via Moore-Penrose generalized inverse of the directed graph's Laplacian. Here in our proposed resistance distance $\Omega$ for Markov chain, it is defined in terms of the group inverse of $I-P$, which according to \cite{M75} is the ``correct generalized inverse to use in connection with finite Markov chains". We leave these open questions above as future work for further comparison between $R$ and $\Omega$.
		\newline
		In \cite[Section $4.2$]{ABPPW15}, the authors propose to view classical effective resistance on \textit{undirected} graph via a variational formula that depends on the modulus. It is unclear to the author whether similar variational formula holds for our proposed effective resistance $\Omega$. The asymmetric nature of $P$ maybe a possible obstacle in generalizing this result to our setting.
		\newline
		Another related work is \cite{BRZ11}. In Section $4$ therein, the authors introduce a few variants of fundamental matrices and express the mean hitting time, commute time as well as the Moore-Penrose inverse of $\Pi(I-P)$ in terms of these fundamental matrices. One can easily express our proposed $\Omega$ in terms of these quantities as well by utilizing Definition \ref{def:omega} and Proposition \ref{prop:equivalentresis}.
		\newline
		Finally, we mention the work \cite{CA02}. In Section $9$ therein, the authors obtain a few interesting relationship between the fundamental matrix $F$ and the group inverse $D$. They can be applied to gain additional insights on these quantities.
	\end{rk}

	We defer the proof of this Proposition to Section \ref{subsec:proofer}. We proceed to investigate whether $\Omega_{i,j}$ defines a metric on $\mathcal{X}$. Recall that in the graph setting its resistance distance $\Omega^G_{i,j}$ defines a metric as its Laplacian $L$ is symmetric. For a proof of this fact one can consult \cite[Section $9.1$]{B10}. This resembles the setting when $P$ is doubly stochastic in which $\Omega_{i,j}$ defines a metric on $\mathcal{X}$, as we shall see in the next Proposition. In general however, $\Omega$ is a semi-metric since it does not satisfy the triangle inequality.
	
	\begin{proposition}\label{prop:rmetric}
		The resistance distance $\Omega_{i,j}$ of $X$, as defined in Definition \ref{def:omega}, satisfies, for $i,j,k \in \mathcal{X}$,
		\begin{enumerate}
			\item(non-negativity)\label{it:nn} $\Omega_{i,j} \geq 0$ and equality holds if and only if $i = j$.
			
			\item(symmetry)\label{it:symmetry} $\Omega_{i,j} = \Omega_{j,i}.$
			
			\item(triangle inequality)\label{it:tri} When $P$ is doubly stochastic, then $\Omega_{i,j} \leq \Omega_{i,k} + \Omega_{k,j}.$
		\end{enumerate}
		In other words, $\Omega = (\Omega_{i,j})_{i,j \in \mathcal{X}}$ defines a semi-metric on $\mathcal{X}$ in general, and is a metric when $P$ is doubly stochastic.
	\end{proposition}

	\begin{rk}[Triangle inequality need not hold for reversible Markov chain]\label{rk:triangle}
		In this Remark, to demonstrate that the triangle inequality (item \eqref{it:tri} in Proposition \eqref{prop:rmetric}) need not hold for reversible finite Markov chain, we provide a simple counterexample by looking at the three-state birth-death Markov chain. Recall that a Markov chain is reversible if and only if it satisfies the detailed balance condition $\pi_i P_{i,j} = \pi_j P_{j,i}$ for all $i,j \in \mathcal{X}$. In this counterexample, suppose that the state space consists of three states with $\mathcal{X} = \{1,2,3\}$, and we consider an ergodic birth-death Markov chain $X$ on $\mathcal{X}$ with birth probability $P_{i,i+1} > 0$ for $i = 1,2$ and death probability $P_{j,j-1} > 0$ for $j = 2,3$. Note that $P_{1,3} = P_{3,1} = 0$ as $X$ is a birth-death chain. It is well-known that birth-death chain is reversible. Now, using Proposition \ref{prop:equivalentresis} we compute
		\begin{align}
			\Omega_{1,3} = \pi_3 \mathbb{E}_1(\tau_3) + \pi_1 \mathbb{E}_3(\tau_1) &= \pi_3 \mathbb{E}_1(\tau_2) + \pi_3 \mathbb{E}_2(\tau_3) + \pi_1 \mathbb{E}_3(\tau_2) + \pi_1 \mathbb{E}_2(\tau_1), \label{eq:omega13}\\
			\Omega_{1,2} + \Omega_{2,3} &= \pi_2 \mathbb{E}_1(\tau_2) + \pi_1 \mathbb{E}_2(\tau_1) + \pi_3 \mathbb{E}_2(\tau_3) + \pi_2 \mathbb{E}_3(\tau_2), \label{eq:omega1223} \\
			\Omega_{1,3} - \Omega_{1,2} - \Omega_{2,3} &= (\pi_3 - \pi_2)\mathbb{E}_1(\tau_2) + (\pi_1 - \pi_2) \mathbb{E}_3(\tau_2), \label{eq:omegadiff}
		\end{align}		
		where we utilize the birth-death property in the second equality of \eqref{eq:omega13}, and \eqref{eq:omegadiff} follows from \eqref{eq:omega13} and \eqref{eq:omega1223}. For three-state birth-death chain with $\pi_1 > \pi_2$ and $\pi_3 > \pi_2$, by \eqref{eq:omegadiff} we then have $\Omega_{1,3} > \Omega_{1,2} + \Omega_{2,3}.$ A concrete numerical example is the following birth-death chain
		$$P = \begin{pmatrix} 
		0.9 & 0.1 & 0 \\ 
		0.5 & 0 & 0.5 \\ 
		0 & 0.1 & 0.9  
		\end{pmatrix}.$$
		Clearly, $P$ is ergodic with $\pi_1 = \pi_3 = 5/11 > 1/11 = \pi_2$. Moreover, we check that $\Omega_{1,3} = 20 > 140/11 = \Omega_{1,2} + \Omega_{2,3}$.
	\end{rk}
	
	The proof of the above Proposition can be found in Section \ref{subsec:proofmetric}. In the following, we present a generalized sum rule of $\Omega$ as one of our major results of this article. The crux of the proof relies on the sum rule of hitting time of Markov chains \cite{PR10} and is deferred to Section \ref{subsec:sumrule}.
	
	\begin{lemma}\label{thm:sumrule}
		Given an ergodic Markov chain $X$ with fundamental matrix $F$ on $\mathcal{X}$, for any square matrices $M, K$ on $\mathcal{X}$ such that
		\begin{enumerate}
			\item $K\1_{|\mathcal{X}|} = \1_{|\mathcal{X}|}$, where $\1_{|\mathcal{X}|}$ is the all-ones vector of length $|\mathcal{X}|$,
			\item $M(K-I)$ is symmetric,
		\end{enumerate} 
		then we have
		$$\sum_{i,j} (M(K-I))_{i,j} \Omega_{i,j} = 2 \mathrm{Tr}(M(I-K)F),$$
		where $\mathrm{Tr}(\cdot)$ is the trace operation.
	\end{lemma}
	
	At first glance, this theorem may seem to be restrictive due to the assumptions on the row sum of $K$ as well as the symmetry of $M(K-I)$. Nonetheless, in many cases these assumptions are fulfilled and we apply the above Lemma \ref{thm:sumrule} which yields the following Corollary on the Markov chain counterpart of Kirchhoff indices:
	
	\begin{corollary}\label{cor:kirchhoff}
	For a given ergodic Markov chain $X$ with non-unit eigenvalues of $P$ given by $(\lambda_i)_{i=2}^{|\mathcal{X}|}$, we have
	\begin{enumerate}
		\item(Kirchhoff index)\label{it:1} $$\sum_{i,j} \Omega_{i,j} = 2 |\mathcal{X}| t^{av} = 2 |\mathcal{X}| \sum_{i=2}^{|\mathcal{X}|} \dfrac{1}{1-\lambda_i}.$$
		
		\item(Multiplicative Kirchhoff index)\label{it:multi} Writing $M$ to be the diagonal matrix with $M_{i,i} = \pi_i$ for all $i \in \mathcal{X}$, 
		$$\sum_{i,j} \pi_i\pi_j\Omega_{i,j} = 2 \mathrm{Tr}(MF-M\Pi).$$
		
		\item(Additive Kirchhoff index)\label{it:add} 
		$$2 t^{av}\leq \sum_{i,j} (\pi_i+\pi_j)\Omega_{i,j} \leq 2 t^{av}(|\mathcal{X}|+1).$$
	\end{enumerate}
	\end{corollary} 

	The above formulae of Markov chain Kirchhoff indices share a striking similarity with their counterparts on graph. For instance, writing $(\lambda_i^L)^{|V|}_{i=2}$ to be the non-zero eigenvalues of the Laplacian $L$, the graph counterpart of Kirchhoff index (see e.g. \cite[Corollary $2$]{PR10}) can be calculated as 
	$$\sum_{i,j} \Omega_{i,j}^G = 2|V|\sum_{i=2}^{|V|} \dfrac{1}{\lambda_i^L},$$
	which resembles the corresponding formula in Corollary \ref{cor:kirchhoff} item \eqref{it:1}.	For recent progress in the study of Kirchhoff indices on graph, we refer interested readers to \cite{P16}. As our second application of the main result of Lemma \ref{thm:sumrule}, we establish a Markov chain counterpart of Foster's first formula of electrical network under a doubly stochastic setting:
	
	\begin{corollary}\label{cor:foster}
		Suppose that $X$ is a reversible Markov chain. By writing $M$ to be the diagonal matrix with $M_{i,i} = \pi_i$ for all $i \in \mathcal{X}$, we then have, for $m \in \mathbb{N}$,
		$$\sum_{i,j} \pi_j P_{j,i}^m \Omega_{i,j} = 2 \mathrm{Tr}\left(M\left(\sum_{j=0}^{m-1}(P^j - \Pi)\right)\right).$$		
		In particular, when $P$ is doubly stochastic, the above gives a Markov chain analogue of the Foster's first formula:
		$$\sum_{i,j} P_{i,j} \Omega_{i,j} = 2(|\mathcal{X}|-1).$$
	\end{corollary}
	
	The above result can be compared to its classical counterpart result in graph theory (see e.g. \cite{PR10}), which gives
	$$\sum_{(i,j)\in E} \Omega^G_{i,j} = |V|-1.$$
	In this vein, we mention the work of \cite{T94} who also gives related results in the direction of Foster's network theorem and reversible Markov chains. 
	
	The rest of the paper is devoted to the proof of the main results. We prove Proposition \ref{prop:equivalentresis} in Section \ref{subsec:proofer}, Proposition \ref{prop:rmetric} in Section \ref{subsec:proofmetric}, Lemma \ref{thm:sumrule} in Section \ref{subsec:sumrule}, Corollary \ref{cor:kirchhoff} in Section \ref{subsec:proofkirchhoff} and finally Corollary \ref{cor:foster} in Section \ref{subsec:foster}.
	
	\section{Proofs of the main results}
	
	\subsection{Proof of Proposition \ref{prop:equivalentresis}}\label{subsec:proofer}
	We first prove item \eqref{it:gir}. It is well-known that $F, \Pi, D$ are connected by the formula $F = \Pi + D$, see e.g. \cite[Theorem $3.1$]{M75}. Desired result follows since
	$$\Omega_{i,j} = F_{i,i} + F_{j,j} - F_{i,j} - F_{j,i} = \pi_i + D_{i,i} + \pi_j + D_{j,j} - \pi_j - D_{i,j} - \pi_i - D_{j,i} = D_{i,i} + D_{j,j} - D_{i,j} - D_{j,i}.$$
	Next, we prove item \eqref{it:mhtr}. Using the relationship
	$\E_i(\tau_j) = \frac{F_{j,j} - F_{i,j}}{\pi_j}$ gives
	$$\Omega_{i,j} = F_{j,j} - F_{i,j} + F_{i,i} - F_{j,i} = \pi_j \mathbb{E}_i(\tau_j) + \pi_i \mathbb{E}_j(\tau_i).$$
	For item \eqref{it:forest}, we only prove the case of $i \neq j$ as the case of $i = j$ is trivial. Using item \eqref{it:mhtr} and \eqref{eq:tree}, we see that
	$$\Omega_{i,j} = \dfrac{q_j}{q} \dfrac{f_{i,j}}{q_j} + \dfrac{q_i}{q} \dfrac{f_{j,i}}{q_i} = \dfrac{f_{i,j} + f_{j,i}}{q}.$$
	Finally, we prove item \eqref{it:ds}. In the doubly stochastic case, $\pi_i = 1/|\mathcal{X}|$ for all $i$, and so by item \eqref{it:mhtr} we write
	$$\Omega_{i,j} = \dfrac{1}{|\mathcal{X}|} \E_i(\tau_j) + \dfrac{1}{|\mathcal{X}|} \E_j(\tau_i) = \dfrac{1}{|\mathcal{X}|}t^{c}_{i,j}.$$
	
	\subsection{Proof of Proposition \ref{prop:rmetric}}\label{subsec:proofmetric}
	We first prove item \eqref{it:nn}. According to Proposition \ref{prop:equivalentresis} item \eqref{it:mhtr}, 
	$$\Omega_{i,j} = \pi_j \mathbb{E}_i(\tau_j) + \pi_i \mathbb{E}_j(\tau_i) \geq 0.$$
	Equality holds if and only if $\E_i(\tau_j) = \E_j(\tau_i) = 0$ if and only if $i = j$. Next, we prove item \eqref{it:symmetry}. Using Proposition \ref{prop:equivalentresis} item \eqref{it:mhtr} again, we have 
	$$\Omega_{i,j} = \pi_j \mathbb{E}_i(\tau_j) + \pi_i \mathbb{E}_j(\tau_i) = \pi_i \mathbb{E}_j(\tau_i) + \pi_j \mathbb{E}_i(\tau_j) = \Omega_{j,i}.$$
	Finally, we prove item \eqref{it:tri} under doubly stochastic $P$. By Proposition \ref{prop:equivalentresis} item \eqref{it:ds}, we see that
	$$\Omega_{i,j} = \dfrac{1}{|\mathcal{X}|}t^{c}_{i,j} \leq \dfrac{1}{|\mathcal{X}|}(t^{c}_{i,k} + t^{c}_{k,j}) = \Omega_{i,k} + \Omega_{k,j},$$
	where we use the triangle inequality for commute time, see e.g. \cite[Chapter $2$, Lemma $9$]{AF14}.
	
	\subsection{Proof of Lemma \ref{thm:sumrule}}\label{subsec:sumrule}
	Using Proposition \ref{prop:equivalentresis} item \eqref{it:mhtr}, we write 
	\begin{align*}
		\sum_{i,j} (M(K-I))_{i,j} \Omega_{i,j} &= \sum_{i,j} (M(K-I))_{j,i} \pi_j \E_i(\tau_j) + \sum_{j,i} (M(K-I))_{i,j} \pi_i \E_j(\tau_i) \\
		&= 2 \mathrm{Tr}(M(I-K)F),
	\end{align*}
	where we use the symmetry of $M(K-I)$ in the first equality, and the second equality follows from the sum rule of the hitting time of Markov chains \cite[Proposition $2$]{PR10}.
	
	\subsection{Proof of Corollary \ref{cor:kirchhoff}}\label{subsec:proofkirchhoff}
	We first prove item \eqref{it:1}. It follows from the random target lemma (see e.g. \cite[Lemma $10.1$]{LPW09}) and Proposition \ref{prop:equivalentresis} item \eqref{it:mhtr} that
	$$\sum_{i,j} \Omega_{i,j} = \sum_{i} \sum_j \pi_j \E_i(\tau_j) + \sum_j \sum_i \pi_i \E_j(\tau_i) = \sum_i t^{av} + \sum_j t^{av} = 2|\mathcal{X}|t^{av} = 2 |\mathcal{X}| \sum_{i=2}^{|\mathcal{X}|} \dfrac{1}{1-\lambda_i},$$
	where the last equality follows from eigentime identity of ergodic Markov chain \cite{CuiMao10}. We proceed to prove item \eqref{it:multi}. In Lemma \ref{thm:sumrule}, by taking $M$ to be the diagonal matrix of the row vector $\pi$ and $K = \Pi$, we readily check that $K\1_{|\mathcal{X}|} = \1_{|\mathcal{X}|}$ and $(M(K-I))_{i,j} = \pi_i \pi_j = (M(K-I))_{j,i}$, and so Lemma \ref{thm:sumrule} gives
	$$\sum_{i,j} \pi_i \pi_j \Omega_{i,j} = 2 \mathrm{Tr}(MF - MKF) = 2 \mathrm{Tr}(MF - M\Pi),$$
	where we use $\Pi F = \Pi$ in the last equality. Finally, we prove item \eqref{it:add}. For the lower bound, applying Proposition \ref{prop:equivalentresis} item \eqref{it:mhtr} again we see that
	$$\sum_{i,j} (\pi_i+\pi_j)\Omega_{i,j} = \sum_{i,j} (\pi_i+\pi_j)(\pi_j \mathbb{E}_i(\tau_j) + \pi_i \mathbb{E}_j(\tau_i)) \geq \sum_{i,j} \pi_i \pi_j \E_i(\tau_j) + \pi_j \pi_i \E_j(\tau_i) = 2t^{av}.$$
	On the other hand, for the upper bound, we have
	\begin{align*}
		\sum_{i,j} (\pi_i+\pi_j)\Omega_{i,j} &= \sum_{i,j} \pi_i^2 \E_j(\tau_i) +  \pi_j^2 \E_i(\tau_j) + \pi_i \pi_j \E_i(\tau_j) + \pi_j \pi_i \E_j(\tau_i) \\
		&\leq \left(\sum_{i,j} \pi_i \E_j(\tau_i) +  \pi_j \E_i(\tau_j)\right) + 2t^{av} = 2|\mathcal{X}|t^{av} + 2t^{av},
	\end{align*}
	where the last equality follows again from the random target lemma.
	
	\subsection{Proof of Corollary \ref{cor:foster}}\label{subsec:foster}
	We first consider 
	\begin{align}\label{eq:fosterpf}
	P^m F = P^m (\Pi + D) = \Pi + \sum_{n =0}^{\infty} (P^{m+n} - \Pi) = \Pi + D - \sum_{j=0}^{m-1}(P^j - \Pi) = F - \sum_{j=0}^{m-1}(P^j - \Pi),
	\end{align}
	where we use $F = \Pi + D$ in the first and last equality. Writing $M$ to be the diagonal matrix of $\pi$ and $K = P^m$, we check that $K\1_{|\mathcal{X}|} = \1_{|\mathcal{X}|}$ and for $i \neq j$ we use the reversibility assumption on $P$ to note that $(M(K-I))_{i,j} = \pi_i P^m_{i,j} = \pi_j P^m_{j,i} = (M(K-I))_{j,i}.$ By Lemma \ref{thm:sumrule}, we have
	$$\sum_{i,j} \pi_j P_{j,i}^m \Omega_{i,j} = 2 \mathrm{Tr}(MF - MP^m F) = 2 \mathrm{Tr}\left(M\left(\sum_{j=0}^{m-1} (P^j - \Pi)\right)\right),$$
	where the last equality follows from \eqref{eq:fosterpf}. In particular, when $P$ is doubly stochastic (and reversible by assumption), its stationary distribution is given by the discrete uniform. As a result, we take $m = 1$ and $\pi_i = 1/|\mathcal{X}|$ to see
	$$\dfrac{1}{|\mathcal{X}|}\sum_{i,j} P_{i,j} \Omega_{i,j} = 2 \mathrm{Tr}(M - M\Pi) = 2\left(1 - \sum_i \pi_i^2\right) = 2 \left(1 - \dfrac{1}{|\mathcal{X}|}\right),$$
	from which the desired result follows. 
	
\noindent \textbf{Acknowledgements}.
The author is grateful to the editor and the anonymous referee for constructive comments that improve the presentation of the manuscript. In particular, the author thanks the referee for pointers to relevant literature in Remark \ref{rk:connections}, for raising the connection with forest representation of mean hitting time that leads to Proposition \ref{prop:equivalentresis} item \eqref{it:forest} and for asking whether triangle inequality holds for reversible chain that leads to Remark \ref{rk:triangle}. The author acknowledges the support from the Chinese University of Hong Kong, Shenzhen grant PF01001143.

\bibliographystyle{abbrvnat}
\bibliography{thesis}

\begin{thebibliography}{27}
\providecommand{\natexlab}[1]{#1}
\providecommand{\url}[1]{\texttt{#1}}
\expandafter\ifx\csname urlstyle\endcsname\relax
  \providecommand{\doi}[1]{doi: #1}\else
  \providecommand{\doi}{doi: \begingroup \urlstyle{rm}\Url}\fi

\bibitem[Albin et~al.(2015)Albin, Brunner, Perez, Poggi-Corradini, and
  Wiens]{ABPPW15}
N.~Albin, M.~Brunner, R.~Perez, P.~Poggi-Corradini, and N.~Wiens.
\newblock Modulus on graphs as a generalization of standard graph theoretic
  quantities.
\newblock \emph{Conform. Geom. Dyn.}, 19:\penalty0 298--317, 2015.

\bibitem[Aldous and Fill(2002)]{AF14}
D.~Aldous and J.~A. Fill.
\newblock Reversible {M}arkov {C}hains and {R}andom {W}alks on {G}raphs, 2002.
\newblock Unfinished monograph, recompiled 2014, available at
  \url{http://www.stat.berkeley.edu/~aldous/RWG/book.html}.

\bibitem[Anantharam and Tsoucas(1989)]{AT89}
V.~Anantharam and P.~Tsoucas.
\newblock A proof of the {M}arkov chain tree theorem.
\newblock \emph{Statist. Probab. Lett.}, 8\penalty0 (2):\penalty0 189--192,
  1989.

\bibitem[Bapat(2010)]{B10}
R.~B. Bapat.
\newblock \emph{Graphs and matrices}.
\newblock Universitext. Springer, London; Hindustan Book Agency, New Delhi,
  2010.

\bibitem[Boley et~al.(2011)Boley, Ranjan, and Zhang]{BRZ11}
D.~Boley, G.~Ranjan, and Z.-L. Zhang.
\newblock Commute times for a directed graph using an asymmetric {L}aplacian.
\newblock \emph{Linear Algebra Appl.}, 435\penalty0 (2):\penalty0 224--242,
  2011.

\bibitem[Chebotarev(2007)]{C07}
P.~Chebotarev.
\newblock A graph theoretic interpretation of the mean first passage times.
\newblock \emph{arXiv preprint math/0701359}, 2007.

\bibitem[Chebotarev and Agaev(2002)]{CA02}
P.~Chebotarev and R.~Agaev.
\newblock Forest matrices around the {L}aplacian matrix.
\newblock \emph{Linear Algebra Appl.}, 356:\penalty0 253--274, 2002.
\newblock Special issue on algebraic graph theory (Edinburgh, 2001).

\bibitem[Chebotarev and Deza(2018)]{CD18}
P.~Chebotarev and E.~Deza.
\newblock Hitting time quasi-metric and its forest representation.
\newblock \emph{Optimization Letters}, Aug 2018.

\bibitem[Coolen-Schrijner and van Doorn(2002)]{CvD02}
P.~Coolen-Schrijner and E.~A. van Doorn.
\newblock The deviation matrix of a continuous-time {M}arkov chain.
\newblock \emph{Probab. Engrg. Inform. Sci.}, 16\penalty0 (3):\penalty0
  351--366, 2002.

\bibitem[Cui and Mao(2010)]{CuiMao10}
H.~Cui and Y.-H. Mao.
\newblock Eigentime identity for asymmetric finite {M}arkov chains.
\newblock \emph{Front. Math. China}, 5\penalty0 (4):\penalty0 623--634, 2010.

\bibitem[Erd\'{e}lyi(1967)]{Erdelyi67}
I.~Erd\'{e}lyi.
\newblock On the matrix equation {$Ax=\lambda Bx$}.
\newblock \emph{J. Math. Anal. Appl.}, 17:\penalty0 119--132, 1967.

\bibitem[Kemeny and Snell(1976)]{KS76}
J.~G. Kemeny and J.~L. Snell.
\newblock \emph{Finite {M}arkov chains}.
\newblock Springer-Verlag, New York-Heidelberg, 1976.
\newblock Reprinting of the 1960 original, Undergraduate Texts in Mathematics.

\bibitem[Kirkland(2014)]{Kirk14}
S.~Kirkland.
\newblock On the {K}emeny constant and stationary distribution vector for a
  {M}arkov chain.
\newblock \emph{Electron. J. Linear Algebra}, 27:\penalty0 354--372, 2014.

\bibitem[Klein(2002)]{K02}
D.~J. Klein.
\newblock Resistance-distance sum rules.
\newblock \emph{Croatica chemica acta}, 75\penalty0 (2):\penalty0 633--649,
  2002.

\bibitem[Levene and Loizou(2002)]{LL02}
M.~Levene and G.~Loizou.
\newblock Kemeny's constant and the random surfer.
\newblock \emph{Amer. Math. Monthly}, 109\penalty0 (8):\penalty0 741--745,
  2002.

\bibitem[Levin et~al.(2009)Levin, Peres, and Wilmer]{LPW09}
D.~A. Levin, Y.~Peres, and E.~L. Wilmer.
\newblock \emph{Markov chains and mixing times}.
\newblock American Mathematical Society, Providence, RI, 2009.

\bibitem[Mao(2004)]{Mao04}
Y.-H. Mao.
\newblock The eigentime identity for continuous-time ergodic {M}arkov chains.
\newblock \emph{J. Appl. Probab.}, 41\penalty0 (4):\penalty0 1071--1080, 2004.

\bibitem[Meyer(1975)]{M75}
C.~D. Meyer, Jr.
\newblock The role of the group generalized inverse in the theory of finite
  {M}arkov chains.
\newblock \emph{SIAM Rev.}, 17:\penalty0 443--464, 1975.

\bibitem[Miclo(2016)]{Miclo16}
L.~Miclo.
\newblock On ergodic diffusions on continuous graphs whose centered resolvent
  admits a trace.
\newblock \emph{J. Math. Anal. Appl.}, 437\penalty0 (2):\penalty0 737--753,
  2016.

\bibitem[Palacios(2016)]{P16}
J.~L. Palacios.
\newblock Some more interplay of the three {K}irchhoffian indices.
\newblock \emph{Linear Algebra Appl.}, 511:\penalty0 421--429, 2016.

\bibitem[Palacios and Renom(2010)]{PR10}
J.~L. Palacios and J.~M. Renom.
\newblock Sum rules for hitting times of {M}arkov chains.
\newblock \emph{Linear Algebra Appl.}, 433\penalty0 (2):\penalty0 491--497,
  2010.

\bibitem[Pitman and Tang(2018)]{PT18}
J.~Pitman and W.~Tang.
\newblock Tree formulas, mean first passage times and {K}emeny's constant of a
  {M}arkov chain.
\newblock \emph{Bernoulli}, 24\penalty0 (3):\penalty0 1942--1972, 2018.

\bibitem[Syski(1978)]{Syski78}
R.~Syski.
\newblock Ergodic potential.
\newblock \emph{Stochastic Process. Appl.}, 7\penalty0 (3):\penalty0 311--336,
  1978.

\bibitem[Tetali(1991)]{Tetali91}
P.~Tetali.
\newblock Random walks and the effective resistance of networks.
\newblock \emph{J. Theoret. Probab.}, 4\penalty0 (1):\penalty0 101--109, 1991.

\bibitem[Tetali(1994)]{T94}
P.~Tetali.
\newblock An extension of {F}oster's network theorem.
\newblock \emph{Combin. Probab. Comput.}, 3\penalty0 (3):\penalty0 421--427,
  1994.

\bibitem[Young et~al.(2016{\natexlab{a}})Young, Scardovi, and Leonard]{YSL16a}
G.~F. Young, L.~Scardovi, and N.~E. Leonard.
\newblock A new notion of effective resistance for directed graphs---{P}art
  {I}: definition and properties.
\newblock \emph{IEEE Trans. Automat. Control}, 61\penalty0 (7):\penalty0
  1727--1736, 2016{\natexlab{a}}.

\bibitem[Young et~al.(2016{\natexlab{b}})Young, Scardovi, and Leonard]{YSL16b}
G.~F. Young, L.~Scardovi, and N.~E. Leonard.
\newblock A new notion of effective resistance for directed graphs---{P}art
  {II}: {C}omputing resistances.
\newblock \emph{IEEE Trans. Automat. Control}, 61\penalty0 (7):\penalty0
  1737--1752, 2016{\natexlab{b}}.

\end{thebibliography}

\end{document}